\documentclass[11pt]{amsart} 
\usepackage{amsfonts}
\usepackage{amsthm}
\usepackage{amsmath}
\usepackage{amssymb}
\usepackage{latexsym}
\usepackage[all]{xy}

\newcounter{zlist}
\newenvironment{zlist}{\begin{list}{{\rm(\arabic{zlist})}}{
\usecounter{zlist}\leftmargin2.5em\labelwidth2em\labelsep0.5em
\topsep0.6ex\itemsep0.3ex plus0.2ex minus0.3ex
\parsep0.3ex plus0.2ex minus0.1ex}}{\end{list}}

\newcounter{blist}
\newenvironment{blist}{\begin{list}{{\rm (\alph{blist})}}{
\usecounter{blist}\leftmargin2.5em\labelwidth2em\labelsep0.5em
\topsep0.6ex\itemsep0.3ex plus0.2ex minus0.3ex
\parsep0.3ex plus0.2ex minus0.1ex}}{\end{list}}

\newcounter{rlist}
\newenvironment{rlist}{\begin{list}{{\rm(\roman{rlist})}}{
\usecounter{rlist}\leftmargin2.5em\labelwidth2em\labelsep0.5em
\topsep0.6ex\itemsep0.3ex plus0.2ex minus0.3ex
\parsep0.3ex plus0.2ex minus0.1ex}}{\end{list}}

\newtheorem{thm}{}[section]

\newcommand{\Ra}{\Rightarrow}
\newcommand{\LRa}{\Leftrightarrow}

\newcommand{\lra}{\longrightarrow}

\newcommand{\Hom}{{\rm Hom}}
\newcommand{\Mor}{{\rm Mor}}

\newcommand{\End}{{\rm End}}
\newcommand{\Fix}{{\rm Fix}}

\newcommand{\Gen}{{\rm Gen}}
\newcommand{\Cogen}{{\rm Cogen}}

\newcommand{\ot}{\otimes}
\newcommand{\id}{{\rm Id}}
\newcommand{\ve}{\varepsilon}

\newcommand{\A}{\mathbb{A}}
\newcommand{\B}{\mathbb{B}}
\newcommand{\M}{\mathbb{M}} 

\newcommand{\wF}{\overline{F}}
\newcommand{\wG}{\overline{G}}

\newcommand{\tF}{\widetilde{F}}
\newcommand{\tG}{\widetilde{G}}

\newcommand{\bS}{{\bf S}} 
\newcommand{\bT}{{\bf T}}

\def\Label{\label}

\begin{document}

\title[Idempotent monads and {\large $\star$}-functors]{Idempotent monads and {\LARGE $\star$}-functors}
\author[J. Clark and R. Wisbauer]{John Clark, Dunedin, New Zealand \\
  Robert Wisbauer, D\"usseldorf, Germany  
  }
\maketitle

\begin{abstract}
For an associative ring $R$, let $P$ be an $R$-module with $S=\End_R(P)$.
C.\ Menini and A.\ Orsatti posed the question of when the related  
functor $\Hom_R(P,-)$ (with left adjoint $P\ot_S-$)  
induces an equivalence between a subcategory of $_R\M$ closed under 
factor modules and a subcategory of $_S\M$ closed under submodules.
They observed that this is precisely the case if the unit of the 
adjunction is an epimorphism and the counit is a monomorphism. 
A module $P$ inducing these properties is called a $\star$-module.

The purpose of this paper is to consider the corresponding question
for a functor $G:\B\to \A$ between arbitrary categories.
We call $G$ a {\em $\star$-functor} if it has a left adjoint $F:\A\to \B$
such that the unit of the adjunction is an {\em extremal epimorphism}
and the counit is an {\em extremal monomorphism}.
In this case $(F,G)$ is an idempotent pair of functors and induces 
an equivalence between the category $\A_{GF}$ of modules for the monad $GF$
and the category $\B^{FG}$ of comodules for the comonad $FG$. Moreover,
$\B^{FG}=\Fix(FG)$ is closed under factor objects in $\B$,  
$\A_{GF}=\Fix(GF)$ is closed under subobjects in $\A$.
\smallskip

Key Words: idempotent monads and comonads, $\star$-modules, equivalence 
of categories, tilting modules, extremal monomorphisms.

AMS classification: 18C15, 16D90
\end{abstract}

\tableofcontents

\section{Introduction}

Let $R$ and $S$ be associative rings and $_RP_S$ an $(R,S)$-bimodule.
In \cite{MenOrs}, C.\ Menini and A.\ Orsatti
asked under which conditions on $P$,
the functors  $P\ot_S-$ and
$\Hom_R(P,-)$ induce an equivalence between 
certain subcategories of $_R\M$ closed under factor modules (i.e.\ $\Gen(P)$) and 
subcategories of $_S\M$ closed under submodules (i.e.\ $\Cogen(\Hom(P,Q))$
for some cogenerator $Q$ in $_R\M$).
Such modules $P$ are called {\em $\star$-modules} and it is well-known that they are closely related to tilting modules (e.g., \cite{ColpiT}, \cite{WiTilt}). 

Because of the effectiveness of these notions in representation theory
of finite dimensional algebras (see Assem \cite{Ass}), various attempts
have been made to extend them to more general situations. This was done mostly in
categories which do permit some technical tools needed (e.g. additivity, tensor product).

The purpose of this article is to filter out the categorical essence of the theory
and to formulate the interesting parts for arbitrary categories. 
For this we consider a pair $(F,G)$ of adjoint functors between categories
$\A$ and $\B$.  
The crucial step is the observation that these induce functors
between the category $\B^{FG}$ of comodules for the comonad
$FG$ on $\B$ and the category $\A_{GF}$ of modules for the monad
$GF$ on $\A$ (see \ref{related.f}). When the comonad $FG$ (equivalently the monad $GF$) is idempotent, 
$\A^{FG}$ may be considered as a coreflective subcategory of $\A$ and
$\B_{GF}$ becomes a reflective subcategory of $\B$ and these categories are 
equivalent. To improve the setting one may additionally require  $\B^{FG}$ to be 
closed under factor objects and $\A_{GF}$ to be closed under subobjects.
This is achieved by stipulating that  
the unit of the adjunction is an {\em extremal epimorphism} in $\A$ and its counit 
is an {\em extremal monomorphism} in $\B$. In this case we say that $G$
is a $\star$-functor or that 
$(F,G)$ is a pair of {\em $\star$-functors}.
Note that no additional structural conditions on the categories are employed. 

By definition, an $(R,S)$-bimodule $P$ is a $\star$-module provided
the functor $\Hom_R(P,-):{_R\M}\to {_S\M}$ is a $\star$-functor and our 
results apply immediately to this situation.  

A $\star$-module $P$ is a {\em tilting module}
if (and only if) $P$ is a subgenerator in $_R\M$. 
To transfer this property to a  $\star$-functor  $G$, one has to 
require that every object $A$ in $\A$ permits a monomorphism $A\to G(B)$
for some $B\in \B$. We will not go into this question here. 

Central to our investigation are the {\em idempotent
monads {\rm (}comonads{\rm )}} which have appeared in various places in the literature,   
e.g., Maranda \cite{Mar}, Isbell \cite{Isb}, Lambek and Rattray 
\cite{LamRat, LamRat.add}, and 
Deleanu, Frei and Hilton \cite{DeHi}.

\section{Preliminaries}

For convenience we recall the basic structures from category theory which will be 
needed in the sequel.

\begin{thm}\Label{monad}{\bf Monads.} \em
A {\em monad} on a category $\A$ is a triple
$\bT=(T,\mu,\eta)$ where $T:\A\to \A$ is an endofunctor and $\mu:TT\to T$,
$\eta:\id_\A \to T$ are natural transformations inducing commutative diagrams
$$\xymatrix{ TTT\ar[r]^{T\mu} \ar[d]_{\mu T} & TT \ar[d]^\mu \\
 TT \ar[r]^\mu & T ,}  \qquad
  \xymatrix{ T\ar[r]^{T\eta} \ar[dr]_{=} & TT \ar[d]^\mu &\ar[dl]^= \ar[l]_{\eta T}  T \\
   &T &. }  $$
\end{thm}

\begin{thm}\Label{mod}{\bf Modules for monads.} \em
Given a monad $\bT=(T,\mu,\eta)$ on the category $\A$, an object $A\in \A$ with a morphism $\rho_A:T(A)\to A$ is called a {\em $\bT$-module} 
(or {\em $\bT$-algebra}) if  
$\rho_A\circ \eta_A = \id_A$ and $\rho_A$ induces commutativity of the diagram
$$\xymatrix{ TT(A)\ar[r]^{T(\rho_A)} \ar[d]_{\mu_A} & T(A) \ar[d]^{\rho_A} \\
 T(A) \ar[r]^{\rho_A}  & A . } $$ 
A {\em morphism} between $\bT$-modules $(A,\rho_A)$ and $(A',\rho_{A'})$ is an  
$f:A\to A'$ in $\A$ 
satisfying $f\circ \rho_A = \rho_{A'} \circ T(f)$.
We denote the set of these morphisms by $\Mor_\bT(A,A')$ and  
the category of $\bT$-modules by $\A_\bT$. 
\end{thm}

\begin{thm}\Label{comonad}{\bf Comonads.} \em
A {\em comonad} on a category $\A$ is a triple
$\bS=(S,\delta,\ve)$ where $S:\A\to \A$ is an endofunctor and $\delta:S\to SS$,
$\ve:S\to \id_\A $ are natural transformations
inducing commutative diagrams
$$\xymatrix{ S\ar[r]^{\delta} \ar[d]_{\delta} & SS \ar[d]^{S\delta} \\
 SS \ar[r]^{\delta S} & SSS, }   \qquad
  \xymatrix{   &\ar[ld]_=  S \ar[d]^\delta \ar[rd]^= & \\
  S & \ar[l]^{S\ve}    SS \ar[r]_{\ve S} &   S.} $$
\end{thm}

\begin{thm}\Label{comod}{\bf Comodules for comonads.} \em
Given a comonad $\bS=(S,\delta,\ve)$ on the category $\A$, 
an object $A\in \A$ with a morphism $\rho^A: A\to S(A)$ is an {\em $\bS$-comodule} if  
$\ve_A\circ \rho^A=\id_A$ and $\rho_A$ induces commutativity of the diagram

$$\xymatrix{ A \ar[r]^{\rho^A} \ar[d]_{\rho^A}& S(A)\ar[d]^{\delta_A} \\
S(A)\ar[r]^{S(\rho^A)}   & SS (A). }$$
A {\em morphism} between $\bS$-comodules $(A,\rho^A)$ and $(A',\rho^{A'})$ is an  
$f:A\to A'$ in $\A$ 
satisfying $\rho^{A'}\circ f = S(f)\circ \rho^A$.
We denote the set of these morphisms by $\Mor^\bS(A,A')$ and    
the category of $\bS$-comodules by $\A^\bS$. 
\end{thm}

\begin{thm}\Label{adjoint}{\bf Adjoint functors.} \em
Let $F:\A\to \B$ and $G:\B\to \A$ be (covariant) functors between
any categories $\A$, $\B$.
 The pair $(F,G)$ is called {\em adjoint} (or an {\em adjunction}) and $F$ (resply.\ $G$) is called a {\em left} (resply.\ {\em right}) {\em adjoint} to $G$ (resply.\ $F$)  if the two
equivalent conditions hold:
\begin{blist}
\item there is an isomorphism, natural in $A\in \A$ and $B\in \B$,
  $$\varphi_{A,B}:\Mor_\B (F(A),B) \to \Mor_\A (A,G(B));$$  
\item there are natural transformations $\eta:\id_\A \to GF$ (called the {\em unit} of the adjunction)
  and  $\ve: FG\to \id_\B$ (called the {\em counit} of the adjunction) with commutative diagrams
  (called the {\em triangular identities})
  $$ \xymatrix{ F \ar[r]^{F\eta} \ar[dr]_= & FGF \ar[d]^{\ve F} \\
                & F }, \quad
     \xymatrix{ G \ar[r]^{\eta G} \ar[dr]_= & GFG \ar[d]^{G\ve} \\
                & G.} $$
\end{blist}
In this case we have the following relations:
$$\begin{array}{rrcl}
\varphi:&  F(A)\stackrel{f}\to B  &\longmapsto & 
A\stackrel{\eta_A} \to GF(A)\stackrel{G(f)}\to G(B) , \\[+1mm]
\varphi^{-1}: &  A\stackrel{g}\to G(B)  & \longmapsto &
 F(A)\stackrel{F(g)}\to FG(B)\stackrel{\ve_B}\to B .
\end{array} $$
\end{thm}

\begin{thm}\Label{adjoint.prop}{\bf Properties of adjoint functors.}   
Let $(F,G)$ be as in {\rm \ref{adjoint}}. Then
 \begin{zlist}
 \item \begin{rlist}
 \item $G$ is faithful if and only if $\ve_B$ is an epimorphism for each $B\in \B$.
 \item $G$ is full if and only if $\ve_B$ is a coretraction {\rm (}split monic{\rm )}
       for each $B\in \B$.
 \item $G$ is full and faithful if and only if $\ve$ is an isomorphism.
 \end{rlist}
 \item \begin{rlist}
 \item $F$ is faithful if and only if $\eta_A$ is a monomorphism for each $A\in \A$.
 \item $F$ is full if and only if $\eta_A$ is a retraction {\rm (}split epic{\rm )}
       for each $A\in \A$.
 \item $F$ is full and faithful if and only if $\eta$ is an isomorphism.  
 \end{rlist}
\end{zlist}  
\end{thm}  

\begin{thm}\Label{adjoint.mon}{\bf Adjoint functors and (co)monads.}    
Let $(F,G)$ be as in {\rm \ref{adjoint}}. Then 
\begin{zlist} 
\item 
\begin{rlist}
\item $\bT=(GF, G\ve F, \eta)$ is a monad on $\A$;
\item there is a functor $\wG:\B\to \A_{GF},\; B\mapsto (G(B), G\ve_B )$.
\end{rlist}
\item  
 \begin{rlist}
 \item $\bS = (FG, F\eta G, \ve)$ is a comonad on $\B$;
 \item there is a functor $\wF:\A\to \B^{FG},\; A\mapsto (F(A), F\eta_A)$.
 \end{rlist}       
\end{zlist}
\end{thm}
\begin{proof}
(1.i), (2.i) are well-known properties of adjoint functors.
\smallskip

(1.ii) describes the {\em comparison functor}. To show its properties 
recall that naturality of $\ve$ yields the commutative
diagram (e.g. \cite[Section 3]{BaWe})
$$\xymatrix{FGFG\ar[r]^{\quad\ve FG} \ar[d]_{FG\ve} & FG \ar[d]^\ve \\
FG \ar[r]^\ve & \id .} $$
Action of $G$ from the left and application to $B$ yields the 
commutative diagram 
$$\xymatrix{ GFGFG(B)\ar[r]^{\quad G\ve FG_B} \ar[d]_{GFG\ve_B} & GFG(B) 
\ar[d]^{G\ve_B} \\
GFG(B) \ar[r]^{G\ve_B} & G(B) .} $$
This proves the associativity condition for the $GF$-module $G(B)$. 
Unitality follows from the triangular identities (\ref{adjoint}).
Again by naturality of
$\ve$, for any $f\in \B$, $G(f)$ is a $GF$-module morphism.  

The proof of (2.ii) is dual to that of (1.ii).
\end{proof}

\begin{thm}\Label{free}{\bf Free functor for a monad.} \em
For any monad $\bT=(T,\mu,\eta)$ on $\A$ and object $A\in \A$,
$(T(A),\mu_A)$ is a $\bT$-module, called the {\em free $\bT$-module} on $A$.
This yields the
{\em free functor}
$$\phi_\bT: \A \to \A_\bT,\quad A\mapsto (T(A),\mu_A),$$
which is left adjoint to the forgetful functor 
$U_\bT: \A_\bT \to \A$ 
 by the isomorphism, for $A\in \A$ and $M\in \A_\bT$,
$$\Mor_\bT (T(A), M) \to \Mor_\A(A,U_\bT(M)),\quad f\mapsto f\circ \eta_A.$$
Notice that $U_\bT \phi_\bT = T$ and $U_\bT(M)=M$ on objects $M\in \A_\bT$. 
The unit of this adjunction is 
  $\eta: \id_\A \to T=U_\bT \phi_\bT$, and for the counit 
  $\tilde\ve: \phi_\bT U_\bT \to \id_{\A_\bT} $
  we have $\mu=U_\bT\tilde\ve \phi_\bT$   
(e.g.\ \cite[Theorem 3.2.1]{BaWe}, \cite[Proposition 4.2.2]{Bor.2}).  
\end{thm}

\begin{thm}\Label{co-free}{\bf Free functor for a comonad.} \em
For any comonad $\bS=(S,\delta,\ve)$ on $\A$ and object $A\in \A$, 
$(S(A),\delta_A)$ is an $S$-comodule, called the {\em free $\bS$-comodule} on $A$.
This yields the {\em free functor}
$$\phi^\bS: \A \to \A^\bS,\quad A\mapsto (S(A),\delta_A),$$
which is right adjoint to the forgetful functor 
$U^\bS: \A^\bS \to \A$ by the isomorphism, for $A\in \A$ and $M\in \A^\bS$,
$$\Mor^\bS ( M ,S(A)) \to \Mor_\A(U^\bS(M),A),\quad g\mapsto \ve_A \circ g.$$
Notice that $U^\bS \phi^\bS = S$ and $U^\bS(M)=M$ on objects in $\A^\bS$. 
The counit of this adjunction is 
  $\ve:U^\bS \phi^\bS=S \to \id_{\A}$, and for the unit 
  $\tilde\eta: \id_{\A^\bS} \to \phi^\bS U^\bS$
  we have $\delta =U^\bS\tilde\eta \phi^\bS$.   
\end{thm}

The following observation is the key to our investigation.

\begin{thm}\Label{id.monad}{\bf Idempotent monads.}  
For a monad $\bT=(T,\mu,\eta)$ on a category $\A$, the following are equivalent:
\begin{blist}
\item The forgetful functor $U_\bT: \A_\bT\to \A$ is full and faithful;
\item the counit $\tilde\ve: \phi_\bT U_\bT \to \id_{\A_\bT}$ is an isomorphism; 
\item the product $\mu:TT\to T$ is an isomorphism;
\item for every $\bT$-module $(A, \rho_A)$, $\rho_A:T(A)\to A$ is an isomorphism in $\A$;
\item $T\eta$ (or $\eta T$) is an isomorphism;
\item $T\eta= \eta T$;
\item $T\mu=\mu T$.
\end{blist}
\end{thm}

\begin{proof} 
A proof of the equivalences from (a) to (d) can be found in \cite[Proposition 4.2.3]{Bor.2}.
The remaining equivalences are shown in \cite[Proposition]{Mar}. Their 
proof is based on the diagram
$$\xymatrix{ TT \ar[r]^\mu \ar[d]_{TT\eta} & T \ar[d]^{T\eta} \\
 TTT \ar[r]^{\mu T} & TT } $$
which is commutative by naturality of $\mu$. 

Now, for example, if  $T\mu=\mu T$, then 
$\mu T\circ TT\eta = \mu T\circ T \eta T = TT$
showing that $\mu$ (and $T\eta$) is an isomorphism, that is, (g)$\Ra$(c). 
\end{proof}

We also need the dual version of this theorem:

\begin{thm}\Label{id.comonad}{\bf Idempotent comonads.}  
For a comonad $\bS=(S,\delta,\ve)$ on a category $\A$, the following are equivalent:
\begin{blist}
\item The forgetful functor $U^\bS: \A^\bS\to \A$ is full and faithful;
\item the unit $\tilde\eta :\id_{\A^\bS} \to \phi^\bS U^\bS  $ is an isomorphism; 
\item the coproduct $\delta:S\to SS$ is an isomorphism;
\item for any $\bS$-comodule $(A, \rho^A)$, $\rho^A:A\to S(A)$ is an isomorphism in $\A$;
\item $S\ve$ (or $\ve S$) is an isomorphism;
\item $S\ve=\ve S$;
\item $S\delta=\delta S$.
\end{blist}
\end{thm}

\section{Idempotent pairs of functors}

In this section, 
we consider an adjoint pair of functors $F:\A\to \B$ and  $G:\B\to \A$
with unit 
$\eta:\id_\A \to GF$ and counit $\ve: FG\to \id_\B$.

\begin{thm}\Label{related.f}{\bf Related functors.}
Let $(F,G)$ be as in {\rm \ref{adjoint}}. 

\begin{zlist}
\item  
For the monad $GF$ on $\A$, composing  
 $U_{GF}$ with $\wF$ {\rm (}from \ref{adjoint.mon}{\rm )} yields a functor 
$$\tF=\wF\circ U_{GF}: \A_{GF} \to \B^{FG}.$$

\item 
For the comonad $FG$ on $\B$, composing $U^{FG}$ with 
$\wG$ {\rm (}from \ref{adjoint.mon}{\rm )} yields a functor

$$\tG =\wG\circ U^{FG}: \B^{FG} \to \A_{GF}.$$

\item  These functors lead to the 
 commutative diagram
$$\xymatrix{ \B^{FG}\ar[rr]^{\tG} \ar[dd]_{U^{FG}} & &
     \A_{GF} \ar[rr]^{\tF}\ar[dd]_{U_{GF}} & & \B^{FG} \ar[dd]^{U^{FG}} \\
     & \\
 \B \ar[rr]_{G}\ar[rruu]^{\wG} & & \A \ar[rr]_{F} \ar[rruu]^{\wF}&& \B ,} $$
\end{zlist}
\end{thm}

In general $(\tF,\tG)$ need not be an adjoint pair of functors. As a first
observation in this context we state:

\begin{thm}\Label{prop.adj}{\bf Proposition.}
Consider an adjoint pair $(F,G)$ {\rm (}as in \ref{adjoint}{\rm )}.
\begin{zlist}
\item For $(A,\rho_A)$ in $\A_{GF}$, the following are equivalent:
\begin{blist}
\item $\eta_A: A\to GF(A)$ is a $GF$-module morphism; 
\item $\eta_A: A\to GF(A)$ is an epimorphism (isomorphism);
\item $\rho_A:GF(A)\to A$ is an isomorphism.
\end{blist}

\item For $(B,\rho^B)$ in $\B^{FG}$, the following are equivalent:
\begin{blist}
\item $\ve_B: FG(B)\to B $ is an $FG$-comodule morphism; 
\item $\ve_B: FG(B)\to B$ is a monomorphism (isomorphism);
\item $\rho^B:B\to FG(B)$ is an isomorphism.
\end{blist}
\end{zlist}
\end{thm}
\begin{proof} (1) (b)$\LRa$(c) for isomorphisms is obvious by unitality of $GF$-modules.

(a)$\Ra$(b)  For $(A,\rho)$ in $\A_{GF}$,
the condition in (a) requires commutativity of the diagram
$$\xymatrix{ GF(A) \ar[r]^{GF\eta_A\;\;} \ar[d]_{\rho_A} 
   & GFGF(A) \ar[d]^{G\ve F(A)} \\  
  A \ar[r]^{ \eta_A } & GF(A) .} $$
By the triangular identities (see \ref{adjoint}), 
$G\ve F \circ GF\eta \simeq \id_{GF}$ and hence $\eta_{A}\circ \rho_A \simeq \id_{G(A)}$.
Since $\rho_A \circ \eta_A \simeq \id_A$ (by unitality)
it follows that 
$\eta_{A}$  (and $\rho_A$) is an isomorphism.

(b)$\Ra$(a) Consider the diagram
$$\xymatrix{ A \ar[d]_{\eta_A} \ar[r]^{\eta_A}& GF(A) \ar[d]^{\eta GF_A} \\
GF(A) \ar[r]^{GF(\eta_A)\quad} \ar[d]_{\rho_A} 
   & GFGF(A) \ar[d]^{G\ve F(A)} \\  
  A \ar[r]^{ \eta_A } & GF(A) ,} $$
in which the upper square is commutative by naturality of $\eta$ and the 
outer rectangle is commutative since the composites of the vertical maps yield 
the identity.
If $\eta_A$ is an epimmorphism, the lower square is also commutative 
showing that $\eta_A$ is a $GF$-module morphism.

(2) These assertions are proved in a similar way.
\end{proof}

\begin{thm}\Label{GF.adj}{\bf $(\tF,\tG)$ as an adjoint pair.}
With the notation in \ref{related.f}, the following are equivalent:
\begin{blist}
\item by restriction and corestriction, $\varphi$ {\rm (}see \ref{adjoint}{\rm )} induces an isomorphism 
$$\tilde\varphi:\Mor^{FG} (\tF(A),B) \to \Mor_{GF} (A,\tG(B))\; 
\mbox{ for } A\in \A_{GF}, \; B\in \B^{FG},$$  
 {\rm (}hence $(\tF,\tG)$ is an adjoint pair of functors{\rm )};
\item $\eta G :G\to GFG$ is an isomorphism;  
\item $G\ve F : GFGF\to GF$ is an isomorphism.
\end{blist}
\end{thm}
\begin{proof}
(a)$\Ra$(b) $\eta_A$ is the image of $\id:\tF(A)\to \tF(A)$ under 
$\tilde\varphi$ and hence a $GF$-module morphism.
By \ref{prop.adj}, this implies that $\eta_A$ is an isomorphism for all
$GF$-modules $A$. Since $G(B)$ is a $GF$-module for any $B\in \B$, 
we have $\eta_{G(B)}:G(B)\to GFG(B)$ an isomorphism, that is, 
$\eta G: G\to GFG$ is an isomorphism.

(b)$\Ra$(c) By the triangular identities, (b) implies 
that $G\ve$ and $G\ve F$ are also isomorphisms. 

(c)$\Ra$(a) Unitality and the triangular identities yield the equalities
$$GF(\rho_A)  \circ GF \eta_A  = G\ve F_A \circ GF \eta_A  =
G\ve F_A  \circ \eta GF_A = \id_{GF_A}.$$
Given (c), we conclude from these that $ GF \eta_A = \eta GF_A$
is an isomorphism and thus $GF(\rho_A)=G\ve F_A.$
With this information, the test diagram for $\eta_A$ being a $GF$-module
morphisms (see proof of \ref{prop.adj}(1)) becomes
$$\xymatrix{ GF(A) \ar[r]^{\eta {GF}_A\;} \ar[d]_{\rho_A} 
   & GFGF(A) \ar[d]^{GF(\rho_A)} \\  
  A \ar[r]^{ \eta_A } & GF(A) ,} $$
and this is commutative by naturality of $\eta$. 
Thus we get an isomorphism  
$$\begin{array}{rcl}
\tilde\varphi:\Mor^{FG} (\tF(A),B)& \lra & \Mor_{GF} (A,\tG(B)), \\
        \tF(A)\stackrel{f}\to B &\longmapsto & 
 A \stackrel{\eta_A}\to \tG\tF(A)\stackrel{\tG(f)}\to \tG(B),
\end{array}
$$ 
showing that $(\tF,\tG)$ is an adjoint pair of functors.
\end{proof}

Adjoint pairs with the properties addressed in \ref{GF.adj} are 
well-known in category theory. Combined with \ref{id.monad} and by standard arguments we obtain the following
list of characterisations for them.

\begin{thm}\Label{GF.id}{\bf Idempotent pair of adjoints.} 
For the adjoint pair of functors $(F,G)$ {\rm (}as in {\rm \ref{adjoint}}{\rm )},
the following are equivalent.
\begin{blist}
\item The forgetful functor $U_{GF}:\A_{GF}\to \A$ is full and faithful;
\item the counit $\bar\ve: \phi_{GF} U_{GF}\to \id_{\A_{GF}}$ is an isomorphism;
\item the product $G\ve F: GFGF \to GF$ is an isomorphism;
\item $\ve F:FGF\to F$ is an isomorphism;
\item the forgetful functor $U^{FG}:\B^{FG}\to \B$ is full and faithful;
\item the unit $\bar\eta: \id_{\B^{FG}} \to \phi^{FG} U^{FG}$ is an isomorphism;
\item the coproduct $F\eta G: FG\to FGFG$ is an isomorphism;
\item $\eta G: G\to GFG$ is an isomorphism.
\end{blist}
{\em If these properties hold then $(F,G)$ is called an}  idempotent 
pair {\em of adjoints}.
\end{thm}

\begin{thm}\Label{remark.id}{\bf Remarks.} \em
Most of these properties have been considered somewhere in the literature. 
Perhaps the first hint of idempotent pairs is given in Maranda 
\cite[Proposition]{Mar} under the name {\em idempotent constructions} (1966).
Isbell discussed their role in \cite{Isb} calling them {\em Galois connections}
(1971).
In Lambek and Rattray \cite{LamRat} they are investigated 
in the context of localisation and duality (1975).  
In the same year they were studied in 
Deleanu, Frei and Hilton \cite[Section 2]{DeHi} where it is shown that
their Kleisli categories are isomorphic to the category of 
fractions (of invertible morphisms). Extending these ideas,
{\em idempotent approximations} to any monad are the topic of 
Casacuberta and Frei \cite{CaFr}.
\end{thm}

For the adjoint functor pair $(F,G)$ we use the notation (e.g.\ \cite{LamRat})
$$\begin{array}{rl}
\Fix(GF,\eta)& = \; \{ A\in \A \,|\, \eta_A:A\to GF(A) \mbox{ is an isomorphism}\},
\\[+1mm]
\Fix(FG,\ve)& = \; \{ B\in \B \,|\, \ve_B:FG(B)\to B \mbox{ is an isomorphism}\}.
\end{array}$$
We denote the (isomorphic) closure of the image of
$GF$ in $\A$ and $FG$ in $\B$ by $GF(\A)$ and $FG(\B)$, respectively.

\begin{thm}\Label{self.theorem}{\bf Idempotent pairs and equivalences.}  
Let $(F,G)$ be an idempotent adjoint pair of functors. Then: 
\begin{rlist}
\item $\A_{GF}\simeq\Fix(GF,\eta)=GF(\A)$ is a reflective subcategory $\A$ with reflector $GF$.
\item $\B^{FG}\simeq\Fix(FG,\ve)=FG(\B)$ is a coreflective subcategory of $\B$ with 
     coreflector $FG$.
\item The (restrictions of the) functors $F$, $G$ 
induce  an equivalence 
 $$ F: GF(\A)  \to FG(\B), \qquad  G : FG(\B)\to GF(\A).$$
\item The Kleisli category of $GF$ is isomorphic to the 
    category of fractions $\A[S^{-1}]$ where $S$ is the family of morphisms 
   of $\A$ rendered invertible by $GF$ (or $F$).
\end{rlist}
\end{thm}

\begin{proof} 
(i) and (ii) follow from \ref{GF.id} (g) and (b), respectively.

(iii) The composition $\tF\tG$ is isomorphic to the identity on $\B^{FG}$ and 
$\tG\tF$ is isomorphic to the identity on $\A_{GF}$.
%

(iv) This is shown in \cite[Theorem 2.6]{DeHi}.
\end{proof}

Of course, if $(F,G)$ induces an equivalence between $\A$ and $\B$, then it is an idempotent pair.
More generally, we obtain from \ref{adjoint.prop} that $(F,G)$ is idempotent
provided the functor $F$ or the functor $G$ is full and faithful.

To consider weaker conditions on the unit and counit, 
recall that an epimorphism $e$ in any category $\A$ is called {\em extremal}
or a {\em cover}
if whenever  $ e=m\circ f$ for a monomorphism $m$ then $m$ is an isomorphism.
Such epimorphisms are isomorphisms if and only if they are monomorph.

\begin{thm}\Label{eta.epi}{\bf $\eta_A$ epimorph.}
Let $(F,G)$ be an adjoint pair of functors {\rm (}as in \ref{adjoint}{\rm )}.
\begin{zlist}
\item If $\eta_A: A\to GF(A)$ is epimorph for any $A\in \A$, then
\begin{rlist}
\item $(F,G)$ is idempotent;
\item $GF$ preserves epimorphisms;
\item for any coproduct $\coprod_{i\in I} A_i$ in $\A$, the canonical 
  morphism $$\psi:{\coprod}_{I} GF(A_i) \to GF({\coprod}_{I} A_i)$$
  is an epimorphism. 
\end{rlist}
\item If 
$\eta_A: A\to GF(A)$ is an extremal epimorphism for any $A\in \A$, then
$\Fix(GF,\eta)$ is closed under subobjects in $\A$.  
\end{zlist}
\end{thm}
\begin{proof}
(1) (i) follows by \ref{prop.adj}. 

(ii) For any morphism $f:A\to A'$ in $\A$, we have the commutative diagram
  $$\xymatrix{A\ar[r]^f \ar[d]_{\eta_A} & A' \ar[d]^{\eta_A'} \\
       GF(A) \ar[r]^{GF(f)} & GF(A'). }$$
If $f$ is epimorph, then so is
the composite $\eta_A'\circ f$ and hence $GF(f)$ must also be epimorph. 

(iii) We have the commutative diagram 
$$\xymatrix{\coprod_{i\in I} GF(A_i)\ar[rr]^\psi   
      && GF(\coprod_{i\in I} A_i )  \\
 &  \coprod_{i\in I} A_i \ar[lu] \ar[ru]_{\eta_{\coprod_I A_i}} } $$
where $\eta_{\coprod_I A_i}$ is epimorph and hence so is $\psi$.

(2) In the diagram in the proof of (1)(ii), assume $f$ to be monomorph and 
$\eta_{A'}$ an isomorphism. Then $\eta_A$ is monomorph and  an extremal
epimorphism which implies that it is an isomorphism.  
\end{proof}

A monomorphism $m$ in any category $\B$ is called {\em extremal} if whenever 
$ m=f\circ e$ for an epimorphism  $e$ then $e$ is an isomorphism.
Such monomorphisms are isomorphisms if and only if they
are epimorph. 

\begin{thm}\label{eps.mono}{\bf $\ve_B$ monomorph.}
Let $(F,G)$ be an adjoint pair of functors (as in \ref{adjoint}).
\begin{zlist}
\item
 Assume  
$\ve_B:FG(B)\to B$ to be monomorph for any $B\in \B$. Then:
\begin{rlist}
\item $(F,G)$ is idempotent;
\item $FG$ preserves monomorphisms;
\item for any product $\prod_{i\in I} B_i$ in $\B$, the canonical 
  morphism $$\varphi:FG({\prod}_{I} B_i) \to {\prod}_{I}FG(B_i)$$
  is a monomorphism. 
\end{rlist}
\item  If
$\ve_B:FG(B)\to B$ is an extremal monomorphism for any $B\in \B$,
then $\Fix(FG,\ve)$ is closed under factor objects in $\B$.  
\end{zlist}
\end{thm}
\begin{proof} The proof is dual to that of \ref{eta.epi}:

(1) (i)  follows by \ref{prop.adj}.

(ii) For any morphism $g:B'\to B$ in $\B$, we have the commutative diagram
  $$\xymatrix{FG(B')\ar[r]^{FG(g)} \ar[d]_{\ve_{B'}} & FG(B) \ar[d]^{\ve_B} \\
       B' \ar[r]^{g} & B. }$$
If $g$ is monomorph, then $g\circ \ve_{B'}$ is monomorph and so is $FG(g)$.

(iii) We have the commutative diagram in $\B$,
$$\xymatrix{ 
FG(\prod_{i\in I} B_i)\ar[rr]^\varphi \ar[rd]_{\ve_{\prod_{I} B_i}} 
      && \prod_{I} FG( A_i) \ar[ld] \\
 &  \prod_{i\in I} B_i ,} $$
where $\ve_{\prod_{I} B_i}$ is monomorph and hence so is $\varphi$.

(2)  In the diagram in (ii), we now have $g$ an epimorphism and $\ve_{B'}$ an isomorphism. Thus $\ve_B$ is epimorph and an extremal monomorphism, 
hence an isomorphism. 
\end{proof}

\begin{thm}{\bf Definition.}\label{def-star} \em
An adjoint pair $(F,G)$ of functors with unit $\eta$ and counit $\ve$
is said to be a {\em pair of $\star$-functors} provided 

\begin{tabular}{l}
$\eta_A: A\to GF(A)$ is an extremal epimorphism for all $A\in \A$ and \\
$\ve_B:FG(B)\to B$ is an  extremal monomorphism for all $B\in \B$.
\end{tabular}
\end{thm}

Combining the information from \ref{self.theorem}, \ref{eta.epi} and \ref{eps.mono}, we obtain the following.

\begin{thm}\Label{th.equiv}{\bf Theorem.}
For a  pair of $\star$-functors $(F,G)$, the functors {\rm (}see \ref{related.f}{\rm )}
$$\tF:\A_{GF} \to \B^{FG}, \quad  \tG: \B^{FG}\to \A_{GF}   $$
induce an equivalence 
where $\A_{GF}=\Fix(GF,\eta)$ is a reflective subcategory of $\A$
closed under subobjects in $\A$ and   
$\B^{FG}=\Fix(FG,\ve)$ is a coreflective subcategory of $\B$ 
closed under factor objects in $\B$.
\end{thm}

\section{$\star$-modules}

In this section let $R$, $S$ be rings and $P$ be an $(R,S)$-bimodule.
The latter provides the adjoint pair of functors
$$T_P:= P\ot_S-: {_S\M}\to {_R\M},\quad H_P:=\Hom_R(P,-): {_R\M}\to {_S\M},$$ 
with unit and counit
$$ \eta_X:  X \to H_P T_P(X),\; 
    x\mapsto [p\mapsto p\ot x], \quad \ve_N:T_PH_P(N)\to N,\; p\ot f \mapsto (p)f,$$ 
where $N\in {_R\M}$ and $X\in {_S\M}$.
Associated to this pair of functors we have the  
 monad and comonad 
$$H_PT_P: {_S\M}\to {_S\M}, \quad T_PH_P:{_R\M}\to {_R\M}.$$
It is well-known that in module categories all monomorphism
and all epimorphisms are extremal.

Recall that $N\in {_R\M}$ is said to be {\em $P$-static} if $\ve_N$ is an isomorphism, and $X\in {_S\M}$ is {\em $P$-adstatic} 
if $\eta_X$ is an isomorphism (e.g.\ \cite{WiStat}).

An $R$-module $N$ is called {\em  $P$-presented}  if there exists an 
exact sequence of $R$-modules
$$ P^{(\Lambda')} \to P^{(\Lambda)} \to N \to 0,\quad \Lambda, \Lambda' 
\mbox{ some sets}.$$

Let $Q$ be any injective cogenerator in $_R\M$ and $P^*:= \Hom_R(P,Q)$. 
An $S$-module $X$ is said to be  {\em  $P^*$-copresented}  
if there exists an exact sequence of $S$-modules
$$ 0\to X \to {P^*}^{\Lambda'} \to {P^*}^{\Lambda},\quad \Lambda, \Lambda' 
\mbox{ some sets}.$$

When $S=\End_R(P)$, there are canonical candidates for fixed modules
for $T_PH_P$ and for $H_PT_P$, namely
$$ P\in \Fix(T_PH_P,\ve) \mbox{ and } S,\; P^* \in \Fix(H_PT_P,\eta),$$
and hence the description of the fixed classes can be related to these objects.

\begin{thm}\Label{HP.id} {\bf $(T_P,H_P)$  idempotent.}
The following are equivalent:
\begin{blist}
\item $H_P\ve T_P: H_PT_P H_P T_P \to H_P T_P$
     is an isomorphism;


\item for any $X\in{_S\M}$, $\ve T_P(X):P\ot_S\Hom(P,P\ot_SX)\to P\ot_S X$ is an isomorphism {\rm (}that is, $P\ot_S X$ is $P$-static{\rm )};


\item $T_P\eta H_P : T_P H_P \to T_P H_P T_P H_P$ is an isomorphism;

\item for any $N\in {_R\M}$, $\eta H_P(N): \Hom(P,N)\to \Hom(P,P\ot_S \Hom(P,N))$ 
 is an isomorphism {\rm (}that is, $\Hom_R(P,N)$ is $P$-adstatic{\rm )}.
\end{blist}

If we assume $S=\End_R(P)$, then {\rm (a)-(d)} are also equivalent to:

\begin{blist}
\setcounter{blist}{4}
\item every $P$-presented $R$-module is $P$-static;
\item every $P^*$-copresented module is $P$-adstatic.
\end{blist}
\end{thm}

\begin{proof} The equivalences (a)-(d) follow from \ref{GF.id}.
For the remaining equivalences see, for example, \cite[4.3]{WiStat}.
\end{proof}

\begin{thm}\Label{id.equiv}{\bf Idempotence and equivalence.}  
With the notation above, 
let $(T_P,H_P)$ be an idempotent pair.  Then these functors 
induce an equivalence 
$$ \widetilde{T_P}:{_S\M}_{H_PT_P} \to {_R\M}^{T_PH_P}, \quad  
 \widetilde{H_P} : {_R\M}^{T_PH_P}\to  {_S\M}_{H_PT_P},$$  
where \,
${_R\M}^{T_PH_P}=\Fix(T_PH_P,\ve)$\, is a coreflective subcategory of ${_R\M}$ and
${_S\M}_{H_PT_P}$ $=\Fix(H_PT_P,\eta)$ is a reflective subcategory of ${_S\M}$:

If $S=\End_R(P)$, then ${_R\M}_{T_PH_P}$ is precisely the subcategory of
$P$-presented $R$-modules and ${_S\M}_{H_PT_P}$ the subcategory of
$P^*$-copresented $S$-modules.
\end{thm}
\begin{proof}
The first part is a special case of \ref{self.theorem}. For the final remark
we again refer to \cite[4.3]{WiStat}.
\end{proof}

Note that the corresponding situation in complete and cocomplete 
abelian categories is described in \cite[Theorem 1.6]{CaGoWi}.

Recall that the module $P$ is {\em self-small} if, for any set $\Lambda$,
the canonical map
$$\Hom_R(P,P)^{(\Lambda)}\to \Hom_R(P,P^{(\Lambda)})$$
is an isomorphism, and $P$ is called 
{\em w-$\Sigma$-quasi\-projective}
 if $\Hom_R(P,-)$ respects exactness of sequences 
   $$0 \to K \to P^{(\Lambda)} \to N \to 0,$$
where $K\in \Gen(P)$, $\Lambda$ any set.

The following observations are essentially known from module theory.

\begin{thm}\Label{etaX.epi}{\bf Proposition.}
For an $R$-module $P$ with $S=\End_S(P)$, the following are equivalent:
\begin{blist}
\item $\eta_X: X\to H_PT_P(X)$ is surjective for all $X\in {_S\M}$;
\item $P$ is self-small and w-$\Sigma$-quasiprojective; 
\item $(T_P,H_P)$ is an idempotent functor pair and ${_S\M_{H_PT_P}}$
  is closed under submodules in ${_S\M}$..
\end{blist}
\end{thm}

For the proof we refer to \cite{WiTilt}, \cite{ColpiE}.  
The assertions where shown by Lambek and Rattray 
for a self-small object in a cocomplete additive category (see
\cite[Theorem 4]{LamRat.add}, \cite[Proposition 1]{Lam.Rem}). 

The following corresponds to \cite[4.4]{WiStat}.

\begin{thm}\Label{ve.mono}{\bf Proposition.}
For an $R$-module $P$ with $S=\End_S(P)$ the following are equivalent:
\begin{blist}
\item $\ve_N: T_PH_P(N)\to N$ is monomorph (injective) for all $N\in {_R\M}$;
\item $(T_P,H_P)$ is idempotent and ${_R\M^{T_PH_P}}$ is closed under factor modules in ${_R\M}$. 
\end{blist}
\end{thm}

As suggested in \ref{def-star}, we call $H_P$ a $\star$-functor
provided the unit $\eta_{_S\M}: \id\to H_PT_P$ is an epimorphism and 
the counit $\ve:T_PH_P\to \id_{_R\M}$ is a monomorphism. 
In this case, the module $P$ is called a {\em $\star$-module}
(\cite{MenOrs}, \cite{ColpiT}) 
and we obtain: 

\begin{thm}\Label{star-mod}{\bf Theorem.}
For an $R$-module $P$ with $S=\End_R(P)$, the following are 
equivalent:
\begin{blist}
\item  $P$ is a  $\star$-module;
\item $H_P$ is a $\star$-functor;
\item $(T_P,H_P)$ induces  an equivalence  
 $$ T_P:{_S\M_{H_PT_P}} \to {_R\M^{T_PH_P}}, \quad  
 H_P : {_R\M^{T_PH_P}}\to  {_S\M_{H_PT_P}},$$
where $_R\M^{T_PH_P}$ is closed under factor modules in $_R\M$
and $_S\M_{H_PT_P}$ is closed under submodules in $_S\M$.
\end{blist}
\end{thm}

The equivalence of (a) and (b) is shown in \cite[Theorem 4.1]{ColpiE}
(see also \cite{MenOrs}, \cite{ColpiT}, \cite{Ass}, \cite{WiTilt}). 
For objects in any Grothendieck category they are shown in Colpi 
\cite[Theorem 3.2]{ColpiG}.
\smallskip

{\bf Acknowledgements.} The authors are grateful to Bachuki Mesablishvili for 
valuable advice and for his interest in their work. 
The research on this topic was started during a visit of the second author 
at the Department of Mathematics at University of Otago in Dunedin, New Zealand.
He wants to express his deeply felt thanks for the warm hospitality and the 
financal support by this institution.

\end{document}